\documentclass[10pt]{article}
\usepackage{amssymb}
\usepackage{amsmath}
\usepackage{amsfonts}
\usepackage{mathrsfs}
\usepackage{amsfonts}
\usepackage{amssymb,amsfonts,amsmath, psfrag,eepic,colordvi,epsfig}
\topskip 
\numberwithin{equation}{section}
\newtheorem{theorem}{Theorem}[section]

\newtheorem{lemma}[theorem]{Lemma}

\begin{document}
	\parskip 7pt
	
	\pagenumbering{arabic}
	\def\sof{\hfill\rule{2mm}{2mm}}
	\def\ls{\leq}
	\def\gs{\geq}
	\def\SS{\mathcal S}
	\def\qq{{\bold q}}
	\def\MM{\mathcal M}
	\def\TT{\mathcal T}
	\def\EE{\mathcal E}
	\def\lsp{\mbox{lsp}}
	\def\rsp{\mbox{rsp}}
	\def\pf{\noindent {\it Proof.} }
	\def\mp{\mbox{pyramid}}
	\def\mb{\mbox{block}}
	\def\mc{\mbox{cross}}
	\def\qed{\hfill \rule{4pt}{7pt}}
	\def\pf{\noindent {\it Proof.} }
	\textheight=22cm

	{\Large
		\begin{center}
			Some identities on the second order  mock theta functions
		\end{center}
	}

\begin{center}
	
	Xingyuan Cai$^{1}$, Eric H. Liu$^2 $  and   Olivia X.M. Yao$^{3}$

	$^{1,3}$School of Mathematical Sciences, \\
	Suzhou University of Science and
	Technology, \\
	Suzhou,  215009, Jiangsu Province,
	P. R. China
	
	$^{2}$School of Statistics and
	Information,\\
	Shanghai University of International
	Business and Economics,\\
	Shanghai, 201620, P. R. China
	
	email: $^{1}$cxy6338@163.com, 
	$^{2}$liuhai@suibe.edu.cn,
	$^{3}$yaoxiangmei@163.com
	
\end{center}

	\noindent {\bf Abstract.}
	 Recently, Nath and Das investigated congruence properties
	  for the second order mock theta function $B(q)$. In their paper,
	   they asked for analytic proofs of three identities
	    on the second order 
	      mock theta functions $A(q)$,  $B(q)$ and $\mu_2(q)$. In this paper, we settle Nath and Das' open problem by 
	    using    the  $(p, k)$-parametrization of theta functions and   several 
	     identities  due to Hickerson and Mortenson.
	

	
	\section{Introduction}
	
	\allowdisplaybreaks

	In   his famous deathbed letter to Hardy,
	Ramanujan gave a list of 17
	functions
	which he called ``mock theta functions".
	He    defined each of
	these functions as a $q$-series
	and found that
	these functions have certain asymptotic
	properties as $q$ approaches
	a root of unity,
	which are similar to theta functions,
	but that they are not theta functions.
	His list contains mock theta
	functions of orders 3, 5 and 7.
	Ramanujan
	recorded some identities for
	mock theta functions of orders
	6 and 10 in his lost notebook.
	After their appearance,
	mock theta functions
	have attracted the
	attention of many mathematicians.
	See \cite{Gordon}
	for a summary  on the  classical
	mock theta functions.
	
	The following   functions are three  classical 
	  mock theta functions of order two 
	 \begin{align*}
	 	A (q):&=\sum_{n=0}^\infty
	 	 a(n)q^n=\sum_{n=0}^\infty
	 	 \frac{(-q;q^2)_n}{(q;q^2)_{n+1}^2}q^{(n+1)^2}
	 	  =\sum_{n=0}^\infty
	 	  \frac{(-q^2;q^2)_n}{(q;q^2)_{n+1}}q^{n+1},  \\
	 	  B (q):&=\sum_{n=0}^\infty
	 	  b(n)q^n=\sum_{n=0}^\infty
	 	  \frac{(-q^2;q^2)_n}{(q;q^2)_{n+1}^2}q^{n(n+1)}
	 	  =\sum_{n=0}^\infty
	 	  \frac{(-q;q^2)_n}{(q;q^2)_{n+1}}q^{n}, \\
	 	  \mu_2 (q):&=\sum_{n=0}^\infty
	 	  \mu (n)q^n=\sum_{n=0}^\infty
	 	   \frac{(-1)^n (q;q^2)_n}{
	 	    (-q^2;q^2)_n^2} q^{n^2}.  
	 \end{align*}
	 Here  and  throughout
	this paper, we use the standard
	$q$-series
	notation:
	\[
	(a;q)_\infty  :=\prod_{k=0}^\infty
	(1-aq^k),\qquad
	(a;q)_n  :=\frac{(a;q)_\infty}
	{(aq^n;q)_\infty}
	\]
	and
	\[
	(a_1,a_2,\ldots, a_m;q)_\infty
	:=(a_1;q)_\infty(a_2;q)_\infty
	\cdots (a_m;q)_\infty,
	\]
	where $q$ is a complex
	number with $|q|<1$. In addition, for any positive integer $k$,
	 we define 
	 \[
	 f_k:=(q^k;q^k)_\infty. 
	 \]

In 2012,  using the theory of modular forms, 
Chan and Mao \cite{Chan-Mao} proved that 
\[
	\sum_{n=0}^\infty b(4n+1)q^n =2\frac{f_2^{10}}{f_1^9}
\]
and 
\[
	\sum_{n=0}^\infty b(4n+2)q^n =2\frac{f_2^{2} f_4^4 }{f_1^5}.
\]
In 2019, Mao \cite{Mao-2019}   proved some congruences modulo 4, 5 and 9 for $b(n)$. For example,
he proved that for $n\geq 0$,
\[
b(10n+6)\equiv b(10n+8) \equiv 0 \pmod 5.
\]
 Mao \cite{Mao} also  proved that 
\[
	\sum_{n=0}^\infty b(6n+2)q^n =4\frac{f_2^{10}f_3^2
	}{f_1^{10}f_6}
\]
and 
\[
	\sum_{n=0}^\infty b(6n+4)q^n =9\frac{f_2^{4} f_3^4 f_6
	}{f_1^8}.
\]
In his nice paper \cite{Wang}, Wang  
 proved that $b(n)$ is odd if and only if $n=2k^2+2k$
  for some $k\geq 0$.
Kaur and Rana \cite{Kaur}  proved that 
 for $n\geq 0$, 
 \[
 b (12n+10)
 \equiv 0 \pmod {36}, \qquad b(18n+16)
 \equiv 0 \pmod {72}. 
 \]
Recently, Yao \cite{Yao} proved some congruences for $B(q)$.
 For example, she  proved that if  $p\equiv 1\ ({\rm mod}\ 4)$,
 then for $n\geq 0$ with $p\nmid n$,
 \begin{align*}
 	b\left(4p^{16\alpha+15}n+\frac{ 5p^{16 \alpha+16}-1
 	}{2}\right)\equiv 0
 	\pmod {32}.
 \end{align*}
 Very recently,
  Nath and Das \cite{Nath} proved that for $n\geq 0$,
  \[
  b(6n+3)
  \equiv 0 \pmod 6,
  \quad
  b(36n+22)\equiv
   0 \pmod {36},
   \quad 
   b(12n+9)\equiv 0 \pmod {54}.
  \]
 To prove the above three congruences, Nath and Das \cite{Nath}
  discovered the following identity   with the help of the $q$-series package developed by Garvan \cite{Garvan}:
  	\begin{align}
  	\sum_{n=0}^\infty
  	b(3n)q^n=\frac{f_2^7f_3^2}{f_1^6f_4f_6}. \label{1-10}
  \end{align}
  They also discovered 
   the following two identities
    related to $A(q)$ and $\mu_2(q)$:
    \begin{align}
    	\sum_{n=0}^\infty
    a(3n+1)q^n &=\frac{f_2^4f_3^2f_4}{f_1^5f_6},
     \label{1-11}\\
    	\sum_{n=0}^\infty
    \mu (3n+1)q^n &=-\frac{f_2^7f_{12}^2
    }{f_1f_4^6f_6}. \label{1-12}
    \end{align}
 
In their paper \cite{Nath}, Nath and Das 
 asked for analytic proofs of \eqref{1-10}--\eqref{1-12}.  The aim
  of this paper is to settle Nath and Das' open problem by
   utilizing  the  $(p, k)$-parametrization of theta functions and several 
    identities   given by Hickerson and Mortenson \cite{Hickerson}.

\section{Analytic proof of \eqref{1-10}}    

To present a analytic proof of 
 \eqref{1-10}, we first prove the following lemma. 

\begin{lemma}\label{L-1}
	Define 
	\begin{align}
		F_1(q):=&	\frac{1}{f_1^4}
		\cdot\frac{f_3}{f_1^3}
		\cdot \frac{f_2^5f_{12}^{14}}{
			f_4f_6f_8^5f_{24}^5}-3q\cdot
		\frac{1}{f_1^4}\cdot
		\left(\frac{1}{f_1f_3}\right)^2 \cdot 
		\frac{f_2^3f_4f_6^5f_{12}^8
		}{f_8^6f_{24}^2}	\nonumber\\
		&+3q\cdot
		\frac{1}{f_1^4}\cdot
		\left(\frac{1}{f_1f_3}\right)^2\cdot
		 \frac{
			f_2f_4^4f_6^{11}f_{24}}{f_8^7f_{12}}
		+\frac{q}{4}\cdot
		\frac{1}{f_1^4}\cdot
		\left(\frac{1}{f_3^4}\right)^2
		\cdot 
		\frac{f_4^8f_6^{26}f_{24}^4}{f_2^4f_8^8f_{12}^{13}}
			\nonumber\\
		&+\frac{3q}{2}\cdot \left(\frac{1}{f_1^4}\right)^2
		\cdot
		f_3^4\cdot \frac{f_2^6f_{12}^{11}}{
			f_6^4f_8^6f_{24}^2}
			 -2q \cdot \left(\frac{1}{f_1^4}\right)^2
		\cdot
	     f_3^4\cdot
	     \frac{f_2^8f_{12}^{20}
		}{f_4^3f_6^{10}f_8^5f_{24}^5} 	
		-\frac{3q}{4}\cdot \frac{f_4^7
			f_{24}^3}{f_8^9} 	\nonumber\\
			& -6q^2\cdot
		\frac{f_3}{f_1^3}\cdot \frac{1}{f_1^4}\cdot 
		\frac{f_2^4f_4^2f_6^2f_{12}^5f_{24}}{f_8^7}
	 +q^2\cdot
		\frac{f_3}{f_1^3}\cdot \frac{1}{f_1^4}\cdot 
		\frac{f_2^2f_4^5f_6^8f_{24}^4}{
			f_8^8f_{12}^4}
			\nonumber\\
		& -2q^2\cdot
		\left(\frac{1}{f_1f_3}\right)^5 \cdot
		\frac{f_4^6f_6^{17}f_{24}^4}{
			f_2f_8^8f_{12}^7}
		  -\frac{q^3}{2}\cdot
		\left( \frac{1}{f_1^4}\right)^2\cdot f_3^4 \cdot 
		 \frac{f_2^5f_4^3f_{12}^2f_{24}^4
		}{f_6f_8^8}-\left(
		\frac{f_3}{f_1^3}\right)^2
		\cdot \frac{f_2^7}{
			f_4f_6}. \label{2-1}
	\end{align}
	Then 
	\begin{align}
	F_1 (q)=0.  \label{2-2}
	\end{align}
\end{lemma}

\noindent{\it Proof.}
It follows from \cite[Entry 25, (v) and (vi), p. 40]{Berndt1991} that 
\begin{align}\label{2-3}
\frac{1}{f_1^4}=\frac{f_4^{14}}{f_2^{14}f_8^4}+4q
 \frac{f_4^2f_8^4}{f_2^{10}}
\end{align}
Replacing $q$ by $q^3$ in \eqref{2-3} yields 
\begin{align}\label{2-4}
\frac{1}{f_3^4}=\frac{f_{12}^{14}}{f_6^{14}f_{24}^4}+4q^3
 \frac{f_{12}^2f_{24}^4}{f_6^{10}}.
\end{align}
Replacing $q$ by $-q$ in \eqref{2-3}, we obtain 
\begin{align}\label{2-4-1}
 f_3^4= \frac{f_{12}^{10}}{
 f_6^2f_{24}^4}-4q^3\frac{f_6^2f_{24}^4
  }{f_{12}^2}.
\end{align}
Baruah and  Ojah \cite[Theorem 4.3]{Baruah}
 proved that 
\begin{align}\label{2-5}
 \frac{1}{f_1f_3}=\frac{f_8^2f_{12}^5}{f_2^2f_4f_6^4f_{24}^2}+q
  \frac{f_4^5f_{24}^2}{f_2^4f_6^2f_8^2f_{12}}
 \end{align}
 In addition, Baruah and  Ojah \cite[Theorem 4.17]{Baruah}  proved that 
\begin{align}\label{2-6}
 \frac{f_3}{f_1^3}=\frac{f_4^6f_6^3}{f_2^9f_{12}^2}+3q
  \frac{f_4^2f_6f_{12}^2}{f_2^7}
 \end{align}
Substituting \eqref{2-3}--\eqref{2-6}  into \eqref{2-1} yields 
\begin{align}\label{2-7}
	F_1(q)=qS_1(q^2)+S_0(q^2),
\end{align}
where
\begin{align}
	S_1(q)=&-2q^3 \frac{f_2^{31}f_3^7f_{12}^{14}}{f_1^{21}f_4^{18}f_6^{12}}
	+16q^3\frac{f_2^{19}f_3f_{12}^8}{f_1^{19}f_4^8}
	+4q^3\frac{f_2^{22}f_3^6f_{12}^{12}}{f_1^{18}f_4^{12}f_6^9}
	-20q^2\frac{f_2^{19}f_3^3f_{12}^6}{f_1^{17}f_4^{10}}
	-8q^2\frac{f_2^7f_6^{12}}{f_1^{15}f_3^3}
		\nonumber\\
	&+64q^2\frac{f_2^{13}f_6^{18}}{f_1^{16}f_3^8f_4^5f_{12}}
	+8q^2\frac{f_2^{10}f_3^2f_6^3f_{12}^4}{f_1^{14}f_4^4}
	-48q^2\frac{f_2^{16}f_6^9f_{12}^2}{f_1^{18}f_3^2f_4^6}
	-q\frac{f_2^{31}f_6^{12}}{2f_1^{23}f_3^3f_4^{16}}
	-32q\frac{f_2f_4^3f_6^{30}}{f_1^{12}f_3^{12}f_{12}^9}
		\nonumber\\
	&-18q\frac{f_2^{18}f_3^3f_6^7f_{12}}{f_1^{17}f_4^{11}}
	+3q\frac{f_2^{21}f_3^9f_{12}^4}{f_1^{19}f_4^{12}f_6^2}
	+4q\frac{f_2^{13}f_3^{11}f_{12}^4}{f_1^{17}f_4^4f_6^6}
	-34q\frac{f_2^7f_6^{12}}{f_1^{13}f_3f_4^2f_{12}^2}
	+3q\frac{f_2^{28}f_3^7f_{12}^5}{f_1^{21}f_4^{15}f_6^3}
		\nonumber\\
	&+24q\frac{f_2^4f_4^2f_6^{21}}{f_1^{14}f_3^6f_{12}^6}
	-3q\frac{f_2^{25}f_3f_6^6f_{12}^2}{f_1^{19}f_4^{14}}
	-\frac{3f_2^7f_{12}^3}{4f_4^9}
	-2\frac{f_2^{25}f_6^{30}}{f_1^{20}f_3^{12}f_4^{13}f_{12}^9}
    +3\frac{f_2^{16}f_3^3f_6^9}{f_1^{17}f_4^7f_{12}^3}
		\nonumber\\
	&+\frac{f_2^{22}f_6^{15}}{4f_1^{18}f_3^2f_4^{12}f_{12}^4}
	+\frac{3f_2^{28}f_6^{21}}{2f_1^{22}f_3^6f_4^{14}f_{12}^6}
	+3\frac{f_2^{15}f_6^{16}}{f_1^{16}f_4^9f_{12}^5}
	+4\frac{f_2^7f_3^2f_6^{12}}{f_1^{14}f_4f_{12}^5}
	-3\frac{f_2^{13}f_6^{18}}{f_1^{15}f_3^3f_4^6f_{12}^6}
	-6\frac{f_2^7f_3^3}{f_1^9} \label{2-8}
\end{align}
and
\begin{align} \label{2-9}
	S_0(q)=&16q^4 \frac{f_2^{10}f_3^6f_{12}^{12}}{f_1^{14}f_4^4f_6^9}
	+32q^4\frac{f_2^7f_3f_{12}^8}{f_1^{15}}
	-96q^3\frac{f_2^4f_4^2f_6^9f_{12}^2}{f_1^{14}f_3^2}
	+128q^3\frac{f_2f_4^3f_6^{18}}{f_1^{12}f_3^8f_{12}}
	-10q^3\frac{f_2^{25}f_3^5f_{12}^{10}}{f_1^{19}f_4^{14}f_6^6}
		\nonumber\\
	&+2q^3\frac{f_2^{31}f_3f_{12}^8}{f_1^{23}f_4^{16}}
	-6q^2\frac{f_2^{28}f_6^9f_{12}^2}{f_1^{22}f_3^2f_4^{14}}
	+8q^2\frac{f_2^{25}f_6^{18}}{f_1^{20}f_3^8f_4^{13}f_{12}}
	+12q^2\frac{f_2^{16}f_3^7f_{12}^5}{f_1^{17}f_4^7f_6^3}
	+2q^2\frac{f_2^{22}f_3^2f_6^3f_{12}^4}{f_1^{18}f_4^{12}}
		\nonumber\\
	&+12q^2\frac{f_2^9f_3^9f_{12}^4}{f_1^{15}f_4^4f_6^2}
	-32q^2\frac{f_2^{13}f_3f_6^6f_{12}^2}{f_1^{15}f_4^6}
	-4q^2\frac{f_2^{19}f_6^{12}}{f_1^{19}f_3^3f_4^8}
	-72q^2\frac{f_2^6f_3^3f_6^7f_{12}}{f_1^{13}f_4^3}
	-9q\frac{f_2^3f_3f_6^4}{f_1^7}
		\nonumber\\
	&+12q\frac{f_2^{16}f_6^{21}}{f_1^{18}f_3^6f_4^6f_{12}^6}
	-16q\frac{f_2^{13}f_6^{30}}{f_1^{16}f_3^{12}f_4^5f_{12}^9}
	+12q\frac{f_2^3f_6^{16}}{f_1^{12}f_4f_{12}^5}
	-6q\frac{f_2^{19}f_6^{12}}{f_1^{17}f_3f_4^{10}f_{12}^2}
		\nonumber\\
	&+12q\frac{f_2^4f_3^3f_4f_6^9}{f_1^{13}f_{12}^3}
	+q\frac{f_2^{10}f_6^{15}}{f_1^{14}f_3^2f_4^4f_{12}^4}
	-14q\frac{f_2f_4^2f_6^{18}}{f_1^{11}f_3^3f_{12}^6}
	+q\frac{f_2^{25}f_3^{11}f_{12}^4}{f_1^{21}f_4^{12}f_6^6}
	-\frac{f_2^{11}f_3^5}{f_1^{11}f_6^4}
	+\frac{f_2^{19}f_3^2f_6^{12}}{f_1^{18}f_4^9f_{12}^5}.
\end{align}
In \cite{Alaca},  Alaca and Williams proved that 
\begin{align}
	f_1=& 2^{-\frac{1}{6}}p^{\frac{1}{24}}(1-p)^{\frac{1}{2}}(1+p)^{\frac{1}{6}}(1+2p)^{\frac{1}{8}}(2+p)^{\frac{1}{8}}k^{\frac{1}{2}}q^{-\frac{1}{24}},\label{2-10} \\
	f_2=&2^{-\frac{1}{3}}p^{\frac{1}{12}}(1-p)^{\frac{1}{4}}(1+p)^{\frac{1}{12}}(1+2p)^{\frac{1}{4}}(2+p)^{\frac{1}{4}}k^{\frac{1}{2}}q^{-\frac{1}{12}}, \label{2-11} \\
	f_3=&2^{-\frac{1}{6}}p^{\frac{1}{8}}(1-p)^{\frac{1}{6}}(1+p)^{\frac{1}{2}}(1+2p)^{\frac{1}{24}}(2+p)^{\frac{1}{24}}k^{\frac{1}{2}}q^{-\frac{1}{8}}, \label{2-12}\\
	f_4=&2^{-\frac{2}{3}}p^{\frac{1}{6}}(1-p)^{\frac{1}{8}}(1+p)^{\frac{1}{24}}(1+2p)^{\frac{1}{8}}(2+p)^{\frac{1}{2}}k^{\frac{1}{2}}q^{-\frac{1}{6}}, \label{2-13} \\
	f_6=&2^{-\frac{1}{3}}p^{\frac{1}{4}}(1-p)^{\frac{1}{12}}(1+p)^{\frac{1}{4}}(1+2p)^{\frac{1}{12}}(2+p)^{\frac{1}{12}}k^{\frac{1}{2}}q^{-\frac{1}{4}}, \label{2-14}\\
	f_{12}=&2^{-\frac{2}{3}}p^{\frac{1}{2}}(1-p)^{\frac{1}{24}}(1+p)^{\frac{1}{8}}(1+2p)^{\frac{1}{24}}(2+p)^{\frac{1}{6}}k^{\frac{1}{2}}q^{-\frac{1}{2}}. \label{2-15}
\end{align}
Substituting \eqref{2-10}--\eqref{2-15} into 
 \eqref{2-8} and \eqref{2-9}, we arrive at 
\begin{align}\label{2-16}
S_0(q)=S_1(q)=0.
\end{align}
Identity \eqref{2-2} follows from \eqref{2-7}
 and \eqref{2-16}. 
This completes the proof of Lemma \ref{L-1}.  \qed

Now, we turn to prove \eqref{1-10}. 
	
	In \cite[(5.2)]{Hickerson},  Hickerson  and Mortenson  proved that 
	\begin{align}
	B(q)=\sum_{n=0}^\infty 
	 b(n) q^n=-\frac{1}{q}m(1,q^4,q^3), \label{2-17}
	\end{align}
where $m(x,q,z) $ is the Appell-Lerch sum defined by
\[
m(x,q,z)=\frac{-z}{(z,q/z,q;q)_\infty }
\sum_{r=-\infty}^\infty
\frac{(-1)^r q^{r(r+1)/2}z^r}{1-xzq^r}.
\]
Hickerson and Mortenson  \cite[Theorem 3.3]{Hickerson} also proved the following 
	identity on Appell-Lerch sums:
	for generic $x,z_0,z_1\in \mathbb{C}^{*} $,
	\begin{align}\label{2-18}
		m(x,q,z_1)-m(x,q,z_0) =&
		\frac{z_0 f_1^3
			j(z_1/z_0;q)j(xz_0z_1;q)}{
			j(z_0;q)j(z_1;q)j(xz_0;q)
			j(xz_1;q)},
	\end{align}
	where 
	\[
	j(z;q)=(z,q/z,q;q)_\infty. 
	\]
	Setting $(x,q,z_1,z_0)\rightarrow
	(1,q^4,q^3,-1)$ in \eqref{2-18} yields 
	\begin{align}\label{2-19}
	m(1,q^4,q^3)-m(1,q^4,-1)=-\frac{f_2^6f_4^3}{
		4f_1^4f_8^4}.
	\end{align}
	In \cite[Corollary 3.8]{Hickerson}, Hickerson and Mortenson   proved that 
	\begin{align}\label{2-20}
		m(x,q,-1)=&m(x^3q^3,q^9,-1)-\frac{x}{q}m(x^3,q^9,-1)
		+\frac{x^2}{q^3}m(x^3/q^3,q^9,-1)\nonumber\\
		&\quad +\frac{xf_1
			f_3^2f_6
			f_9 j(x^2q;q^2) }{2qf_2^2
			f_{18}^2
			j(-x^3;q^3)}.
	\end{align}
	Taking $(x,q) \rightarrow (1,q^4)$ in \eqref{2-20}, we arrive at 
	\begin{align}\label{2-21}
	m(1,q^4,-1)=m(q^{12},q^{36},-1)
	-q^{-4}m(1,q^{36},-1)+q^{-12}m(q^{-12},q^{36},-1)
	+\frac{ f_4^3f_{12}^3
		 f_{36}}{4q^4f_8^3f_{24}f_{72}^2}.
	\end{align}
In view of   \eqref{2-17},
	 \eqref{2-19} and \eqref{2-21}, 
	\begin{align}
		\sum_{n=0}^\infty b(n)q^n&=-\frac{1}{q}\left(
		m(1,q^4,-1)  -\frac{f_2^6f_4^3}{
			4f_1^4f_8^4}\right)\nonumber\\
		&=-\frac{1}{q}m(q^{12},q^{36},-1)
		+q^{-5}m(1,q^{36},-1)-q^{-13}m(q^{-12},q^{36},-1)
		\nonumber\\
		&\qquad 		-\frac{q^{-5}}{4}
		\cdot\frac{f_4^2}{f_8}
		\cdot\frac{f_4}{f_8^2}\cdot \frac{ f_{12}^3
			 f_{36}}{f_{24}f_{72}^2}+
		\frac{ q^{-1}}{4}
		\cdot
		\left( \frac{f_2f_4}{f_1f_8}\right)^4
		\cdot \frac{f_2^2}{f_4}.\label{2-22}
	\end{align}
It follows from \cite[Corollary (i), p. 49]{Berndt1991} that 
\begin{align}\label{3-11}
	\frac{f_2^5}{f_1^2f_4^2}=\frac{f_{18}^5}{f_9^2f_{36}^2
	}+2q\frac{f_6^2f_9f_{36}}{f_3f_{12}f_{18}}. 
\end{align}
Replacing $q$ by $-q$ in \eqref{3-11}
yields 
	\begin{align}\label{2-23}
	\frac{f_1^2}{f_2}=\frac{f_9^2}{f_{18}}
	-2q\frac{f_3f_{18}^2}{f_6f_9}.
	\end{align}
	Hirschhorn and Sellers \cite{Hirschhorn} proved that 
	\begin{align}\label{3-12}
		\frac{f_2}{f_1f_4}=\frac{f_{18}^9}{
			f_3^2f_9^3f_{12}^2f_{36}^3}+q\frac{f_6^2f_{18}^3}{
			f_3^3f_{12}^3}+q^2\frac{f_6^4f_9^3f_{36}^3}{
			f_3^4f_{12}^4f_{18}^3}.
	\end{align}
	Replacing $q$ by $-q$ in \eqref{3-12}
	yields 
	\begin{align}\label{2-24}
	\frac{f_1}{f_2^2}=\frac{f_3^2f_9^3}{
		f_6^6}-q\frac{f_3^3f_{18}^3}{
		f_6^7}+q^2\frac{f_3^4f_{18}^6 }{f_6^8f_9^3}.
	\end{align}
	In \cite{Hirschhorn-2014}, Hirschhorn and Sellers  proved the following 
	identity:
	\begin{align}\label{2-25}
	\frac{f_2f_4}{f_1f_8}
	=\frac{f_6^2f_9f_{36}^6}
	{f_3^2f_{12}f_{18}^3f_{24}f_{72}^2}
	+q\frac{f_{12}^2f_{18}^6f_{72}}{
		f_3f_6f_9^2f_{24}^2f_{36}^3}+q^2\frac{
		f_6f_9f_{12}f_{72}}{f_3^2f_{24}^2}.
	\end{align}
	Substituting  \eqref{2-23}, \eqref{2-24}
	 and \eqref{2-25}
	  into \eqref{2-22} and  extracting  those
	  terms in which the power of  $q
	  $ is congruent to 0 modulo 3, then
	  replacing 
	  $q^{3}$ by $q$, we arrive at 
	\begin{align*}
		\sum_{n=0}^\infty
		b(3n )q^n= &	\frac{f_2^5f_3f_{12}^{14}}{
			f_1^7f_4f_6f_8^5f_{24}^5}-3q\frac{f_2^3f_4f_6^5f_{12}^8
		}{f_1^6f_3^2f_8^6f_{24}^2}+3q\frac{
			f_2f_4^4f_6^{11}f_{24}}{f_1^6f_3^2f_8^7f_{12}}
		+\frac{q}{4}\frac{f_4^8f_6^{26}f_{24}^4}{f_1^4f_2^4f_3^8f_8^8f_{12}^{13}}
		\nonumber\\
		& +\frac{3q}{2}\frac{f_2^6f_3^4f_{12}^{11}}{
			f_1^8f_6^4f_8^6f_{24}^2}-2q\frac{f_2^8f_3^4f_{12}^{20}
		}{f_1^8f_4^3f_6^{10}f_8^5f_{24}^5}-\frac{3q}{4}\frac{f_4^7
			f_{24}^3}{f_8^9}-6q^2
		\frac{f_2^4f_3f_4^2f_6^2f_{12}^5f_{24}}{f_1^7f_8^7}
		\nonumber\\
		&+q^2\frac{f_2^2f_3f_4^5f_6^8f_{24}^4}{
			f_1^7f_8^8f_{12}^4}-2q^2
		\frac{f_4^6f_6^{17}f_{24}^4}{
			f_1^5f_2f_3^5f_8^8f_{12}^7}
		-\frac{q^3}{2}\frac{f_2^5f_3^4f_4^3f_{12}^2f_{24}^4
		}{f_1^8f_6f_8^8}. 
	\end{align*}
	In light of \eqref{2-1} and the above identity,
	\begin{align}\label{2-26}
	\sum_{n=0}^\infty
	b(3n )q^n=F_1(q)+\frac{f_2^7f_3^2}{f_1^6f_4f_6}.
	\end{align}
Identity \eqref{1-10} follows from \eqref{2-2} and \eqref{2-26}.This completes the proof. \qed

\section{Analytic proof of \eqref{1-11}}    

	 We first prove the following lemma. 
	 
	  	\begin{lemma}\label{L-2}
	 	Define
	 	\begin{align}
	 		F_2(q):=&	2q \cdot
	 		\frac{1}{f_1^4}\cdot
	 		\frac{f_1}{f_3^3}\cdot \frac{ f_4^2f_6^{11}
	 			f_{24}^2}{f_2^2f_8^2f_{12}^5}+2q\cdot
	 		\frac{f_3^3}{f_1  } \cdot \frac{1}{f_1^4}\cdot \frac{f_2^3f_{12}^7}{
	 			f_4^2f_6^4f_8f_{24}}	\nonumber\\
	 		&+\frac{1}{2}\cdot
	 		\frac{1}{f_1^4}\cdot \frac{f_2^2f_{12}^{10}}{
	 			f_4^3f_6f_{24}^4}
	 		-\frac{f_2^2f_8f_{12}^9}{2f_4^6f_6^3f_{24}^3}+\frac{f_6^3}{f_4^3}
	 		-\left(\frac{f_3}{f_1^3}\right)^2\cdot \frac{f_2^4f_4}{f_6}.
	 		\label{3-1}
	 	\end{align}
	 	Then
	 	\begin{align}
	 		F_2(q)=0. \label{3-2}
	 	\end{align}
	 	
	 \end{lemma}
	  	\noindent{\it Proof.}  
	Xia and Yao \cite{Xia} proved that 
	\begin{align}\label{3-3}
	\frac{f_3^3}{f_1}=\frac{f_4^3f_6^2}{f_2^2f_{12}}+q\frac{ f_{12}^3
	}{ f_4}.
	\end{align}
Replacing $q$ by $-q$ in \eqref{3-3} yields 
	\begin{align}\label{3-4}
	\frac{f_1}{f_3^3}=	\frac{f_2f_4^2f_{12}^2}{f_6^7}-q\frac{f_2^3 f_{12}^6
	}{ f_4^2f_6^9}.
	\end{align}
		Substituting \eqref{2-3}, \eqref{2-6}, \eqref{3-3}
	and \eqref{3-4} into  \eqref{3-1} yields
	\begin{align}
		F_2(q)=q H_1(q^2)+H_0(q^2), \label{3-5}
	\end{align}
	where
	\begin{align*}
		H_1(q)=&-8q \frac{f_2^2f_3^2f_4^2f_6f_{12}^2}{f_1^9}
		+8q\frac{f_4^3f_6^{10}}{f_1^7f_2f_3^4f_{12}}
		+2\frac{f_2^{18}f_3^4f_{12}^2}{f_1^{15}f_4^6f_6^3}+2\frac{f_2^{15}f_6^6}{f_1^{13}f_3^2f_4^5f_{12}}	 	\nonumber\\
		&+2\frac{f_4^4f_6^{10}}{f_1^8f_2f_3f_{12}^4}
		-6 \frac{f_2^9f_3^3}{f_1^{12}}
	\end{align*}
	and
	\begin{align*}
		H_0(q)=&-2q \frac{f_2^{14}f_3^2f_6f_{12}^2}{f_1^{13}f_4^6}
		+8q\frac{f_2^6f_3^4f_4^2f_{12}^2}{f_1^{11}f_6^3}
		+2q\frac{f_2^{11}f_6^{10}}{f_1^{11}f_3^4f_4^5f_{12}}+8q\frac{f_2^3f_4^3f_6^6}{f_1^9f_3^2f_{12}}
		-9q\frac{f_2^5f_3f_6^4}{f_1^{10}}
		\nonumber\\
		&+\frac{f_2^{11}f_6^{10}}{2f_1^{12}f_3f_4^4f_{12}^4}-\frac{f_1^2f_4f_6^9}{2f_2^6f_3^3f_{12}^3}+\frac{f_3^3}{f_2^3}-\frac{f_2^{13}f_3^5}{f_1^{14}f_6^4}.
	\end{align*}
	Substituting \eqref{2-10}--\eqref{2-15}
	into the above two identities,
	we see that 	 
	\begin{align}\label{3-6}
		H_0(q)=H_1(q)=0.	
	\end{align}
	Identity \eqref{3-2} follows from \eqref{3-5} and \eqref{3-6}. This completes the proof of Lemma \ref{L-2}.\qed 	
	 	     
Now, we are ready to prove 
 \eqref{1-11}.

	Taking $(x,q,z_1,z_0) \rightarrow
	(q,q^4,q^2,-1)$
	in \eqref{2-18} leads to  
	\begin{align}\label{3-8}
	m(q,q^4,q^2)-m(q,q^4,-1)
	=-\frac{f_4^9}{2f_1f_2^3f_8^4}.
	\end{align}
	Setting $(x,q)\rightarrow (q,q^4)$
	in \eqref{2-20}, we arrive at 
	\begin{align}\label{3-9}
		m(q,q^4,-1)=m(q^{15},q^{36},-1)
		-q^{-3}m(q^3,q^{36},-1)+q^{-10}m(q^{-9},q^{36},-1)
		+\frac{f_2f_3f_{12}^2f_{24}f_{36}}{2q^3f_6^2f_8f_{72}^2}.
	\end{align}   
	In \cite{Hickerson}, Hickerson and 
	 Mortenson proved that 
	 \[
	 A(q)=\sum_{n=0}^\infty a(n)q^n=-m(q,q^4,q^2),
	 \]
	 from which with
	  \eqref{3-8} and \eqref{3-9},
	   we have 
	\begin{align}\label{3-10}
		\sum_{n=0}^\infty a(n)q^n= &\frac{f_4^9}{2f_1f_2^3f_8^4}-m(q,q^4,-1)
		\nonumber\\
		=& 
		  -m(q^{15},q^{36},-1)
		+q^{-3}m(q^3,q^{36},-1)-q^{-10}m(q^{-9},q^{36},-1)
		\nonumber\\
		&+\frac{1}{2}\cdot
		\frac{f_2}{f_1f_4}\cdot
		\left(\frac{f_4^{5}}{f_2^2f_8^2}\right)^2-\frac{1}{2q^3} \cdot \frac{f_2}{f_4^2} \cdot 
		\frac{f_4^2}{f_8} \cdot \frac{f_3f_{12}^2f_{24}f_{36}}{ f_6^2f_{72}^2}.
	\end{align}
	Substituting 
	  \eqref{3-11}--\eqref{2-24}   into \eqref{3-10}
	  and extracting 
	    those
	  terms in which the power of  $q
	  $ is congruent to 1 modulo 3, then
	  dividing by $q$ and replacing  
	  $q^{3}$ by $q$, we obtain 
		\[
	\sum_{n=0}^\infty
	a(3n+1)q^n =\left(F_2(q)+\frac{f_2^4f_3^2f_4}{f_1^6f_6}\right)f_1,
	\]
where $F_2(q)$ is defined by \eqref{3-1}. 	
 Identity 
  \eqref{1-11}
   follows from \eqref{3-2}
    and the above identity. This completes
     the proof. \qed 
	
\section{Analytic proof of \eqref{1-12}}    

To prove \eqref{1-12}, we first 
 prove the following lemma.

	\begin{lemma}\label{L-4}
		Define  
		\begin{align}
			F_3(q):=&-4q\cdot \frac{f_1^3}{f_3}\cdot \frac{f_6^4f_{12}^2}{f_2^4f_4^2}
			+4\frac{f_2^3f_6^3}{f_4^4}-4\cdot 
			(f_1f_3)^2\cdot \frac{
				f_6^4}{f_2^4f_4f_{12}}+\left(\frac{f_1^2}{f_3^2}\right)^2\cdot
			\frac{f_6^{13}}{
				f_2^7f_{12}^4} 
			\nonumber\\
			&-2\frac{f_2^5f_8f_{12}^9}{f_4^7f_6^3f_{24}^3}+\frac{f_2^7f_{12}^2}{f_4^6f_6}.\label{4-1}
		\end{align}
		Then 
		\begin{align*}
			F_3(q)=0. 
		\end{align*}
	\end{lemma}
	\noindent{\it Proof.} 
	 Xia and Yao  \cite[Lemma 2.4]{Xia} proved  that 
	\begin{align}\label{4-2}
		\frac{f_1^3}{f_3}=\frac{f_4^3}{f_{12}}-3q
		\frac{f_2^2f_{12}^3}{f_4f_6^2}.
	\end{align}
It follows from \cite[(30.10.2) and (30.12.3)]{Hirschhorn-2017}
	  that 
	\begin{align}\label{4-3}
		\frac{f_1^2}{f_3^2}=\frac{f_2f_4^2f_{12}^4}{f_6^5f_8f_{24}}-2q
		\frac{f_2^2f_8f_{12}f_{24}}{f_4f_6^4}
	\end{align}
	and 	\begin{align}\label{4-4}
		f_1f_3=\frac{f_2f_8^2f_{12}^4}{f_4^2f_6f_{24}^2}-q
		\frac{f_4^4f_6f_{24}^2}{f_2f_8^2f_{12}^2}.
	\end{align}
	Substituting \eqref{4-2}--\eqref{4-4}  into \eqref{4-1} yields 
	\begin{align}\label{4-5}
		F_3(q)=R_0(q^2),
	\end{align}
	where
	\begin{align} \label{4-6}
		R_0(q)=&12q \frac{f_3^2f_6^5}{f_1^2f_2^3}
		-4q\frac{f_2^7f_3^6f_{12}^4}{f_1^6f_4^4f_6^5}
		+4q\frac{f_3^5f_4^2f_{12}^2}{f_1^3f_2^2f_6^2}
		+4\frac{f_1^3f_3^3}{f_2^4}
		-4\frac{f_3^2f_4^4f_6^7}{f_1^2f_2^5f_{12}^4}
		\nonumber\\
		&+\frac{f_2^4f_3^3f_6^4}{f_1^5f_4^2f_{12}^2}
		-2\frac{f_1^5f_4f_6^9}{f_2^7f_3^3f_{12}^3}+\frac{f_1^7f_6^2}{f_2^6f_3}.
	\end{align}
	Substituting \eqref{2-10}--\eqref{2-15} into 
	\eqref{4-6}, we arrive at 
	\begin{align}\label{4-7}
		R_0(q)=0.
	\end{align}
Lemma \ref{L-4}  follows from \eqref{4-5}
	and \eqref{4-7}. 
	This completes the proof of Lemma \ref{L-4}.  \qed

   To conclude this section, we present a 
  proof
   of \eqref{1-12}.

In \cite{Hickerson}, Hickerson and Mortenson 
 proved that 
 \begin{align}\label{4-8}
 	\mu_2(q)&=\sum_{n=0}^\infty
 	\mu(n)q^n =4m(-q,q^4,-1)-\frac{f_2^8}{f_1^3f_4^4}.
 \end{align}
 Replacing $q$ by $-q$ in \eqref{3-9} yields
 	\begin{align*} 
 	m(-q,q^4,-1)=&m(-q^{15},q^{36},-1)
 	+q^{-3}m(-q^3,q^{36},-1) 
 		\nonumber\\
 	&+q^{-10}m(-q^{-9},q^{36},-1)-\frac{f_2f_6f_{12}f_{24}f_{36}}{2q^3f_3f_8f_{72}^2},
 \end{align*} 
 from which with \eqref{4-8}, we arrive at 
 \begin{align}\label{4-9}
 	\sum_{n=0}^\infty
 	\mu(n)q^n =& 4m(-q^{15},q^{36},-1)
 	+4q^{-3}m(-q^3,q^{36},-1)+4q^{-10}m(-q^{-9},q^{36},-1)
 	\nonumber\\
 	&-\frac{2}{q^3}
 	\cdot
 	\frac{f_2}{f_4^2}
 	\cdot
 	\frac{f_4^2}{f_8}\cdot \frac{f_6f_{12}f_{24}f_{36}}{f_3f_{72}^2}
 	-\left(\frac{f_2^5}{f_1^2f_4^2}\right)^2
 	\cdot
 	 \frac{f_1}{f_2^2}.
 \end{align}
 Substituting \eqref{3-11}, 
\eqref{2-23} and  \eqref{2-24}
   into  \eqref{4-9}
 and extracting 
 those
 terms in which the power of  $q
 $ is congruent to 1 modulo 3, then
 dividing by $q$  and replacing  
 $q^{3}$ by $q$, we obtain 
 \begin{align}
 \sum_{n=0}^\infty
 \mu(3n+1)q^n =\frac{1}{f_1}\left(F_3(q)-\frac{f_2^7f_{12}^2
 }{f_4^6f_6}\right). \label{4-10}
 \end{align}
 Identity \eqref{1-12}
  follows from 
  Lemma \ref{L-4} and \eqref{4-10}. This completes
   the proof. \qed

  \section*{Statements and Declarations}

 \noindent{\bf Funding.}
 \noindent{\bf Acknowledgments}
This work was supported by 
the National Natural Science Foundation of
China  (grant
12371334).

 \noindent{\bf Competing Interests.}
 The authors declare that they have  
 no conflict of interest.

 \noindent{\bf Data Availability Statements.} Data sharing not applicable to this
 article as no datasets were generated or analysed during the current
 study.


\begin{thebibliography}{HD}
		
 
 \bibitem{Alaca}
 S. Alaca and K.S. Williams, 
 The number of representations of a positive integer by certain octonary quadratic forms, Funct. Approx. Comment. Math.  43 (2010)  45--54. 
 
 
 	\bibitem{Baruah}
 N.D. Baruah
  and  K.K. Ojah, 
		Analogues of Ramanujan's partition identities and 
		congruences arising from his theta functions and modular equations,
		 Ramanujan J.
		  28 (2012) 385--407.
		
		\bibitem{Berndt1991}
		B.C. Berndt, Ramanujan's Notebooks, Part III, Springer, New York,
		1991.
		
		
	 
		
		\bibitem{Chan-Mao}
		S.H. Chan and
		R.R. Mao, Two congruences for
		Appell-Lerch sums, Int.
		J. Number Theory 8  (2012) 111--123.
		
		\bibitem{Garvan}
		F. Garvan,  A $q$-product tutorial for a $q$-series maple package,
		  S\'{e}m. Lothar. Combin. 42 (1999) 27--52.
		  
		 		\bibitem{Gordon} 
		  B. Gordon and R.J. McIntosh, A survey of classical mock theta functions, in
		  Developments in Mathematics. Developments in Mathematics, Volume 23 New York:
		  Springer, pp. 95--144, (2012).
		  
		  
			\bibitem{Hickerson}
		D.R. Hickerson
		 and E.T. Mortenson,    Hecke-type double sums, Appell-Lerch sums,
		  and mock theta functions, I, 
		  Proc. London Math. Soc. 109 (2014) 382--422.
		  
	\bibitem{Hirschhorn-2017}
	  M.D. Hirschhorn,  The power of $q$, a personal journey, Developments in Mathematics, v.
		  49, Springer, 2017.
		  
		  	\bibitem{Hirschhorn}
		 M.D. Hirschhorn and J.A. Sellers,  Arithmetic properties of partitions with odd parts distinct, Ramanujan J. 22 (2010) 273--284.
		 
		 
		  	\bibitem{Hirschhorn-2014}
	M.D. Hirschhorn and J.A. Sellers, A congruence modulo 3 for partitions into distinct non-multiples of four, J. Integer Sequences 17 (2014) Article 14.9.6.
		  
		\bibitem{Kaur} H.
		  Kaur and M. Rana, On second order mock theta function $B(q)$,  Electron. Rese. Arch. 30 (2022) 52--65. 
		  
		  		  
			\bibitem{Mao-2019}
		R.R. Mao, Two identities on the mock theta function $V_0(q)$,
		  J. Math. Anal. Appl. 479 (2019) 122--134. 
		  
		  
		
			\bibitem{Mao}
		R.R. Mao, Arithmetic properties of coefficients of the mock theta function $B(q)$, Bull. Austr. Math. Soc. 102 (2020) 50--58.
 
 	\bibitem{Nath}
 H. Nath and H.  Das, 	Infinite families 
 of congruences 
 for the second order mock theta 
 function $B(q)$, submitted (arXiv: 2509.20708v1).
 
 
 \bibitem{Wang}
 L.Q. Wang,  Parity of coefficients of mock theta functions, J. Number Theory 229  (2021)  53--99.
		
		\bibitem{Xia}
		E.X.W. Xia and O.X.M. Yao, New Ramanujan-like congruences modulo powers of 2 and 3 for
		overpartitions, J. Number Theory 133 (2013) 1932--1949.
		
		
		\bibitem{Yao}
	O.X.M. Yao, Congruences for the second order mock theta function $B(q)$,
	 Rev. Real Acad. Cienc. Exactas F\'{i}s. Nat. Ser. A Mat.   119 (2025) \# 112.
	 
 
	   
	   
	 
		
	\end{thebibliography}
\end{document}